\documentclass[11pt]{amsart2000}
\newtheorem{theorem}{Theorem}
\newtheorem{lemma}[theorem]{Lemma}
\newtheorem{corollary}[theorem]{Corollary}
\newtheorem{example}[theorem]{Example}
\theoremstyle{definition}
\newtheorem{definition}[theorem]{Definition}
\newtheorem{notation}[theorem]{Notation}
\newtheorem*{noproblem}{Problem}

\begin{document}
\setlength{\unitlength}{0.01in}
\linethickness{0.01in}
\begin{center}
\begin{picture}(474,66)(0,0) 
\multiput(0,66)(1,0){40}{\line(0,-1){24}}
\multiput(43,65)(1,-1){24}{\line(0,-1){40}}
\multiput(1,39)(1,-1){40}{\line(1,0){24}}
\multiput(70,2)(1,1){24}{\line(0,1){40}}
\multiput(72,0)(1,1){24}{\line(1,0){40}}
\multiput(97,66)(1,0){40}{\line(0,-1){40}} 
\put(143,66){\makebox(0,0)[tl]{\footnotesize Proceedings of the Ninth Prague Topological Symposium}}
\put(143,50){\makebox(0,0)[tl]{\footnotesize Contributed papers from the symposium held in}}
\put(143,34){\makebox(0,0)[tl]{\footnotesize Prague, Czech Republic, August 19--25, 2001}}
\end{picture}
\end{center}
\vspace{0.25in}
\setcounter{page}{265}
\title{Wallman-Frink proximities}
\author{Marlon C. Rayburn}
\address{Department of Mathematics\\
The University of Manitoba\\
Winnipeg, Manitoba\\
R3T 2N2\\
Canada}
\email{rayburn@cc.umanitoba.ca}
\thanks{Marlon C. Rayburn,
{\em Wallman-Frink proximities},
Proceedings of the Ninth Prague Topological Symposium, (Prague, 2001),
pp.~265--270, Topology Atlas, Toronto, 2002}
\begin{abstract}
This is a survey of compactification extension results and problems for a 
special class of proximities.
\end{abstract}
\maketitle

The compatible Efremovi\v c proximities on a Tichonov space are ordered by
$\delta \leq \rho$ if for all $A, B \subseteq X$, 
$A \delta B \rightarrow A \rho B$.
It is well known that the Smirnov completions of the Efremovi\v c
proximities give a one-to-one, order reversing correspondence between the
Hausdorff compactifications and the proximities. 
This merely means that the proximities share the compactification lattice
order problem (``find a necessary and sufficient condition on $X$ that
the Hausdorff compactifications form a lattice'').

While it is not practical to work on the order problems by embedding the 
Hausdorff compactifications in the larger family of $T_1$ compactifications,
there simply being no end to the latter, the Efremovi\v c proximities
generalize nicely to the complete lattice of Ladato proximities. 
When Gagrat and Naimpally showed that the compatible separated proximities
on a space complete to $T_1$ compactifications, it seemed a solution to
the order problem was at last in sight.

Unfortunately, troubles remain. 
The Gagrat-Naimpally compactifications are among the $T_1$
compactifications, yes, but which ones are they?
Moreover, the corresponding order between the completions breaks down. 
Even a $p$-map between proximity spaces will only lift, in general, to a
continuous extension from the Gagrat-Naimpally completion of the domain to
the ``bunch space'' of the range --- something rich and strange.

What follows is an attempt to solve these problems for the special case of
Lodato proximities in the style of Wallman-Frink.

\begin{definition}
Let $X$ be a $T_1$ space with at least two points.
Let $\mathcal{B}$ be any base for the closed sets such that 
\begin{enumerate}
\item
$\mathcal{B}$ is a network (i.e. $x \in G$,
open, implies there is some
$B \in {\mathcal{B}}$ with $x \in B \subseteq G)$, and
\item
$\mathcal{B}$ is a ring of sets (i.e. $\mathcal{B}$ is
closed under finite unions and
finite intersections).
\end{enumerate}
Then $\mathcal{B}$ is a {\it normal} base for the closed sets.
\end{definition}

N.~b.\ since $x$ has at least two points, 
$\emptyset \in {\mathcal{B}}$.

\begin{definition}
The {\it Wallman-Frink} proximity $\delta_b$ associated with normal base
$\mathcal{B}$ is given by: 
$$A \not\!\!{\delta}_b B 
\equiv 
\exists F_1, F_2 \in {\mathcal{B}}
\mbox{ with } A \subseteq F_1, B \subseteq F_2 
\mbox{ and }
F_1 \cap F_2 = \emptyset.$$

Hence 
$A \delta_b B \iff 
\forall F_1, F_2 \in {\mathcal{B}}, A \subseteq F_1 
\mbox{ and }
B \subseteq F_2 \Rightarrow F_1 \cap F_2 \not = \emptyset$.
\end{definition}

\begin{theorem}
$\delta_b$ is a compatible, separated Lodato proximity on $X$.
\end{theorem}

\begin{proof}
Since $\emptyset \in {\mathcal{B}}$, clearly 
$\emptyset \not\!\!{\delta}_b A$ 
for any $A \subseteq X$.
It is also clear that 
$A \not\!\!{\delta}_b B \Rightarrow B \not\!\!{\delta}_b A$,
and 
$A \not\!\!{\delta}_b B \Rightarrow A \cap B = \emptyset$.

Since $\mathcal{B}$ is a ring of sets, if 
$A \not\!\!{\delta}_b B$ and 
$A \not\!\!{\delta}_b C$, then 
$A \not\!\!{\delta}_b (B \cup C)$.

Since $\mathcal{B}$ is a network, 
$x \not = y \Rightarrow x \not\!\!{\delta}_b Y$,
and $x \not \in cl(A) \iff x \not\!\!{\delta}_b A$. 
Thus 
$A \not\!\!{\delta}_b C$ and 
$B \subseteq cl(C) \Rightarrow \not\!\!{\delta}_b B$.
\end{proof}

\begin{corollary}
$\delta_b$ is Efremovi\v c $\iff$ whenever $F_1$, $F_2 \in {\mathcal{B}}$
such that $F_1 \cap F_2 = \emptyset$, then there exists some 
$C$, $D \in {\mathcal{B}}$ such that $F_1 \subseteq X \backslash C$,
$F_2 \subseteq X \backslash D$, and 
$(X \backslash C) \cap (X \backslash D) = \emptyset$.
\end{corollary}

A {\it b-filter} is a filterbase $\mathcal{F}$ of sets from $\mathcal{B}$
such that whenever
$B \in {\mathcal{F}}$ and $B \subseteq F \in {\mathcal{B}}$, then $F \in
{\mathcal{F}}$. A
{\it b-ultrafilter} is a maximal b-filter. By Zorn's Lemma, every
b-filter is contained in at least one
b-ultrafilter. N.b., for any ultrafilter $\mu$ on $X$, $\mu
\cap {\mathcal{B}}$ is a b-filter.

\begin{definition}
An ultrafilter $\mu$ on $X$ is an {\it ultrafilter of type b}
if $\mu \cap {\mathcal{B}}$ is a b-ultrafilter.
\end{definition}

\begin{lemma}
For any {\it normal} base $\mathcal{B}$ on $T_1$-space $X$ and each
$x \in X$, the point ultrafilter $\mu_x = \{A\colon x \in A\}$ is
an ultrafilter of type $b$.
\end{lemma}

\begin{proof} 
We need only check that the b-filter $\mu_x \cap {\mathcal{B}}$ is a
maximal among the b-filters. 
Suppose $\varphi$ is a b-filter with 
$\mu_x \cap {\mathcal{B}} \subseteq \varphi$. 
Then each 
$F \in \varphi \subseteq {\mathcal{B}}$, 
$F \cap B \not = \emptyset$ for all
$B \in \mu_x \cap {\mathcal{B}}$. 
Since $\mathcal{B}$ is a network, we see that $x \in cl(F) = F$ for each
$F \in \varphi$. 
Hence $\mu_x \cap {\mathcal{B}} =\varphi$.
\end{proof}

\begin{lemma}
If $\mu$ is any b-ultrafilter, then there is an ultrafilter $\alpha$ of
type $b$ such that $\mu = \alpha \cap {\mathcal{B}}$.
\end{lemma}

\begin{proof} 
$\mu$ is a filterbase of sets, so it is contained an ultrafilter $\alpha$. 
Thus $\mu \subseteq \alpha \cap {\mathcal{B}}$. 
Since $\mu$ is maximal, $\mu =\alpha \cap {\mathcal{B}}$.
\end{proof}

\begin{lemma}
Let $\delta$ be any compatible Lodato proximity with $\delta_b \leq \delta$.
For every ultrafilter $\mu$ on $X$, there is an ultrafilter $\alpha$ of
type $b$ such that $\mu \delta \alpha$.
\end{lemma}

\begin{proof} 
Since $\mu \cap {\mathcal{B}}$ is a b-filter, we must have 
$\mu \cap {\mathcal{B}} \subseteq \varphi$ for some b-ultrafilter 
$\varphi$. 
By the last lemma, there is an ultrafilter $\alpha$ of type $b$ with
$\varphi = \alpha \cap {\mathcal{B}}$. 
Suppose $\mu \not\!{\delta} \alpha$.
Then $\mu \not\!\!{\delta}_b \alpha$, so we must have an 
$M \in \mu$, and $A \in \alpha$ and some $F_1$, $F_2 \in {\mathcal{B}}$
with $M \subseteq F_1$, $A \subseteq F_2$ and $F_1 \cap F_2 = \emptyset$. 
But since $\mu \cap {\mathcal{B}} \subseteq \alpha \cap {\mathcal{B}}$, we
have $F_1$ and $F_2 \subseteq \alpha$,
contradiction to filter.
\end{proof}

\begin{lemma}
On $(X, \delta_b)$, given any ultrafilter $\mu$ and any ultrafilter 
$\alpha$ of type $b$, $\mu \delta_b \alpha$ if and only if 
$\mu \cap {\mathcal{B}} \subseteq \alpha \cap {\mathcal{B}}$.
Hence any two ultrafilters near an ultrafilter of type $b$ are near each
other. 
In particular, $\delta_b$ is transitive over ultrafilters of type $b$.
\end{lemma}

\begin{proof}
Suppose $\mu \delta_b \alpha$. 
Then for each $B \in \mu \cap {\mathcal{B}}$, $B \cap F \not= \emptyset$
or every $F \in \alpha \cap {\mathcal{B}}$.
Since $\alpha \cap {\mathcal{B}}$ is a b-ultrafilter, we have that 
$B \in \alpha \cap {\mathcal{B}}$. 
The converse is clear from the definition of $\delta_b$.
\end{proof}

\begin{definition}
Let $\gamma$ be any filterbase on proximity space $(X, \delta)$ and set 
$\Pi_\delta(\gamma)=\bigcup\{\upsilon \colon \upsilon \hbox{ is an 
ultrafilter and } \upsilon \delta \gamma\}$.
\end{definition}

\begin{definition}
\mbox{}
\begin{itemize}
\item[a)]
A grill $\gamma$ on a proximity space $(X, \delta)$ is a {\it precluster}
if whenever $A \subseteq X$ and $\{A\} \delta \gamma$, then 
$A \in \gamma$.
\item[b)]
A pre-cluster $\sigma$ on $(X,\delta)$ is a {\it cluster} if it is a clan.
\end{itemize}
\end{definition}

\begin{example}
Let $\alpha$ be any filterbase on proximity space $(X, \delta)$. 
Then $\Pi(\varphi)$ is a pre-cluster.
\end{example}

\begin{theorem}
On $(X, \delta_b)$, if $\alpha$ is an ultrafilter of type $b$, then
$\Pi(\alpha)$ is a cluster.
\end{theorem}

\begin{proof}
By Example 12, $\Pi(\alpha)$ is a pre-cluster and by lemma 9,
$\Pi(\alpha)$ is a clan.
\end{proof}

\begin{theorem}
The subspace $T_bX$ of the Gagrat-Naimpally completion 
$\alpha_{\delta_b}X$ given by the set of all maximal clans is a $T_1$
compactification on $X$.
\end{theorem}

\begin{proof}
By Theorem 11, 
$$T_bX = \{\Pi(\alpha) \colon 
\mbox{$\alpha$ is an ultrafilter of type $b$}\}$$ 
is a subset of $\alpha_\delta X$, the set of all maximal clans
on $X$. 
Hence $T_bX$ is a $T_1$ space, which by Lemma 6 contains a dense copy of
$X$. It remains to show that $T_bX$ is compact.

For each $F \in {\mathcal{B}}$, let 
$$F'=\{\Pi(\alpha)\colon F \in \Pi(\alpha)\}.$$
These will be the basic closed sets of the topology of $T_bX$. 
Let $\mathcal{L} = \{F'_j \colon j \in \Gamma\}$ be a family of basic
closed sets with the finite intersection property. 
Let 
$${\mathcal{F}} = \{F_j \colon j \in \Gamma\} \subseteq {\mathcal{B}}.$$
Let 
$${\mathcal{F}}' = 
\{ \bigcap_{j \in \forall} F_j \colon \Lambda \hbox{ a non-empty, finite 
subset of } \Gamma \}.$$
Then ${\mathcal{F}}'$ is a filterbase of sets from ${\mathcal{B}}$ since 
$\mathcal{B}$ is a ring, and ${\mathcal{F}} \subseteq {\mathcal{F}}'$. 
By Zorn's Lemma, ${\mathcal{F}}'$ is contained in some b-ultrafilter,
which by Lemma 7 we may write as $\alpha \cap {\mathcal{B}}$ of some
ultrafilter $\alpha$ of type $b$. 
But then ${\mathcal{F}} \subseteq \Pi(\alpha)$, so for each 
$j \in \Gamma$, $\Pi(\alpha) \in F'_j$, and thus 
$\bigcap\mathcal{L} \not = \emptyset$.
\end{proof}

If ${\mathcal{B}}$ is a {\it normal} base of closed sets for $T_1$ space
$X$, we may constuct a ``Wallman-Frink'' compactification $bX$ of $X$ in
the usual way (See Willard's ``General Topology'', Exercise 19K p 142):
\begin{itemize}
\item[a)] 
Let $bX$ be the set of all b-ultrafilters on $X$.
\item[b)] 
For each $B \in {\mathcal{B}}$, let 
$B' = \{\varphi \in bX \colon B \in \varphi\}$.
This is a base for closed sets for a topology on $bX$. 
Call it the ``absorption'' topology for $\mathcal{B}$.
\item[c)] 
Under the absorption topology, $bX$ is compact and $T_1$. 
\item[d)] 
$\tau \colon X \to bX$ by $t(x) = \mu_x \cap \mathcal{B}$ is a 
homeomorphism from $X$ onto a dense subset of $bX$.
\end{itemize}

\begin{theorem}
There is a homeomorphism between the $T_1$ compactification $T_bX$ of 
$(X,\delta_b)$ and the Wallman-Frink compactification $bX$, which fixes
$X$ point-wise.
\end{theorem}

\begin{proof}
Let $f$ be the bijection assigning maximal b-clan $\sigma = \Pi(\alpha)$
to the b-ultrafilter $\sigma \cap {\mathcal{B}}$. 
Set $B \in {\mathcal{B}}$ gives us both basic closed set
$$
B* = \{\sigma \in T_bX \colon B \in \sigma \}
\ \mbox{and} \
B' = \{\varphi
\in bX \colon B \in \varphi\}.
$$
But clearly $B \in \varphi$ if and only if $B \delta \alpha$ for every
ultrafilter $\alpha$ of type $b$ for which 
$\varphi = \alpha \cap {\mathcal{B}}$, and $F \in \alpha$ if and only if
$F \in \sigma_\alpha$ for 
$\sigma_\alpha = 
\bigcup\{\mu \colon \hbox{ ultrafilter } \mu \delta_b \alpha \}$. 
Hence $f$ is one-to-one between basic closed sets of $bX$ and $T_bX$,
and is therefore a homeomorphism. That $f$ is $1-1$ on $X$ follows from
Lemma 6.
\end{proof}

\begin{noproblem}
This identifies the Wallman compactification of a $T_1$ space with a
subspace of the Gagrat-Naimpally compactification of $(X, \delta_w)$. 
Conceivably this is possible, yet it does not rule out the equality of the
two extensions. 
This would amount to showing the reverse of Theorem 11, that every maximal
clan is of form $\Pi(\alpha)$ for some ultrafilter $\alpha$ of type $b$.
\end{noproblem}

Specifically, let $\sigma$ be any maximal clan on $(X, \delta_b)$
and consider 
$$
\sigma \cap {\mathcal{B}} = 
\bigcap \{\mu \cap {\mathcal{B}}\colon \hbox{ ultrafilter } \mu \subseteq
\sigma\}.
$$
Let $F \in \sigma \cap {\mathcal{B}}$ and $F \subseteq B \in
{\mathcal{B}}$. 
Then there is some ultrafilter $\mu \subseteq \sigma$ such that 
$F \in \mu$. 
Therefore $B \in \mu$, so $B \in \sigma \cap {\mathcal{B}}$. 
Now let $F_1$ and $F_2 \in \sigma \cap {\mathcal{B}}$.
Then $F_1 \delta_b F_2$ since $\sigma$ is a clan, so 
$F_1 \cap F_2 \not = \emptyset$.
Suppose $F_1 \cap F_2 \in \sigma$. 
Then $\sigma \cap {\mathcal{B}}$ is a b-filter, so there is an 
ultra-filter $\alpha$ of type $b$ such that 
$\sigma \cap \mathcal{B} \subseteq \alpha \cap {\mathcal{B}}$.
Hence for every $\mu \subseteq \sigma$, $\mu \delta_b \alpha$, so
$\sigma \subseteq \Pi(\alpha)$. 
By the maximality of $\sigma$, $\sigma = \Pi(\alpha)$. 
Thus

\begin{theorem}
Let $\sigma$ be a maximal clan on $(X, \delta_b)$. 
Then $\sigma = \Pi(\alpha)$ for some ultrafilter $\alpha$ of type $b$ 
$\iff$ 
$[F_1, F_2 \in \sigma \Rightarrow F_1 \cap F_2 \in \sigma]$.
\end{theorem}

Now suppose $X$ to be a $T_1$-space for which $\mathcal{B}$ and 
$\mathcal{D}$ are {\it normal} bases with 
${\mathcal{D}} \subseteq {\mathcal{B}}$. 
[Hence $\delta_b \leq \delta_d$.]

\begin{lemma}
Let $\sigma$ be any element of $T_d$. Then
there is an ultrafilter
$\beta$ of both types $b$ and $d$ such that $\sigma = \Pi_d(\beta)$.
\end{lemma}

\begin{proof}
Since $\sigma$ is an element of $T_d$, there is an ultrafilter $\gamma$ of 
type $d$ such that $\sigma = \Pi_d(\gamma)$. 
But $\gamma$ is an ultrafilter in $(X, \delta_b)$, so there exists an
ultrafilter $\beta$ of type $b$ such that $\gamma \delta_b \beta$.
This means that $\gamma \cap {\mathcal{B}} \subseteq \beta$. 
But $\mathcal{D} \subseteq {\mathcal{B}}$ implies
$\gamma \cap {\mathcal{D}} \subseteq {\mathcal{B}} \cap {\mathcal{D}}$. 
By maximality of $\gamma \cap {\mathcal{D}}$, we have 
$\gamma \cap {\mathcal{D}} = \beta \cap {\mathcal{D}}$. 
Thus $\beta$ is also an ultrafilter of type $d$.

Now ultrafilter $\mu \subseteq \Pi_d(\gamma)$ if and only if 
$\mu \cap {\mathcal{D}} \subseteq 
\gamma \cap {\mathcal{D}} = 
\beta \cap {\mathcal{D}}$. 
Therefore $\mu \subseteq \Pi_d(\gamma)$ if and only if 
$\mu \subseteq \Pi_d(\beta)$, and $\sigma = \Pi_d(\beta)$.
\end{proof}

\begin{notation}
An ultrafilter $\beta$ of both types $b$ and $d$ will be called an
ultrafilter of joint type.
\end{notation}

\begin{lemma}
Let $\alpha$ be any ultrafilter of type $b$.
Then there is an ultrafilter $\beta$ of joint type such that $\alpha
\delta_d \beta$.
\end{lemma}

\begin{proof}
There is an ultrafilter $\gamma$ of type $d$ such that 
$\sigma \delta_d \gamma$. 
By the last lemma, there is an ultrafilter $\beta$ of joint type with 
$\gamma \delta_d \beta$.
Hence $\alpha \delta_d \beta$.
\end{proof}

\begin{definition}
Let $\beta$ be an ultrafilter of joint type and put 
$$\Gamma_\beta = \{\Pi_b(\alpha) \colon 
\mbox{$\alpha$ is an ultrafilter of type $b$ and $\alpha \delta_d \beta$} 
\}.$$
\end{definition}

\begin{theorem}
$\Gamma_\beta$ is a $\delta_d$-clan, and is contained in a unique maximal
$\delta_d$-clan, namely $\Pi_d(\beta)$.
\end{theorem}

\begin{proof}
First, $\Gamma_\beta$ is a $d$-clan. 
Let $S$ and $T$ be sets in $\Gamma_\beta$.
By definition of $\Gamma_\beta$, we can find ultrafilters $\alpha_1$ and
$\alpha_2$ so that $\alpha_1 \delta_d \beta$, $\alpha_2 \delta_d \beta$
and $S \in \alpha_1$, $T \in \alpha_2$. 
Suppose $S \not\!\!{\delta}_d T$. 
Then there exist $D_1$, $D_2 \in {\mathcal{D}}$ such that 
$S \subseteq D_1$, $T \subseteq D_2$ and $D_1 \cap D_2 = \emptyset$. 
Since $\alpha_1$ and $\alpha_2$ are ultrafilters, we must have 
$D_1 \in \alpha_1$ and $D_2 \in \alpha_2$.
Now $\alpha_1 \cap {\mathcal{D}} \subseteq \beta \cap {\mathcal{D}}$ and
$\alpha_2 \cap {\mathcal{D}} \subseteq \beta \cap {\mathcal{D}}$. 
Thus $D_1$ and $D_2$ are in $\beta \cap {\mathcal{D}}$, contradiction to
$D_1 \cap D_2 = \emptyset$.
Clearly $\Gamma_\beta \subseteq \Pi_d(\beta)$, which is a maximal
$d$-clan. 
Suppose $\Pi$ is a maximal $d$-clan with 
$\Gamma_\beta \subseteq \Pi \not = \Pi_d(\beta)$. 
Then there must an $A \in \Pi$ and some $G \in \beta$ for which 
$A \not\!\!{\delta}_d G$. 
But $A \delta_d B$ for all $B \in \Gamma_\beta$, and, because 
$\beta \subseteq \Pi_b(\beta) \subseteq \Gamma_\beta$,
$G \in \Gamma$, contradiction. 
Thus $\Pi_d(\gamma)$ is unique.
\end{proof}

\begin{theorem}
The following are clearly equivalent:
\begin{itemize}
\item[a)] 
For each $\beta$ of joint type, $\Gamma_\beta = \Pi_d(\beta)$. \item[b)]
If $\alpha$ is of type $b$ and $\beta_1$, $\beta_2$ are of joint type such
that $\beta_1 \delta_d \alpha$ and $\beta_2 \delta_d \alpha$, then
$\beta_1 \delta_d \beta_2$.
\item[c)]
$\{\Gamma_\beta \colon \beta \hbox{ is of joint type}\}$ partitions
$T_bX$.
\end{itemize}
\end{theorem}

\begin{theorem}
If any (hence, all) of the conditions of $22$ are met, then there exists a
continuous function $f \colon T_bX \to T_dX$ which is the identity on $X$.
\end{theorem}

\begin{proof}
Define $f \colon T_bX \to T_dX$ by the following:
For $\Pi_b(\alpha) \in T_bX$, there is some $\beta$ of joint type such
that $\alpha \delta_d\beta$.
Let $f[\Pi_b(\alpha)]=\Pi_d(\beta)$. 
By 19, $f$ is a well defined function.

To show that $f$ is continuous, let ${\mathcal{A}} \subseteq T_bX$ and
suppose $\sigma_0 \in T_bX$ with 
$f(\sigma_0) \not \in cl_df[{\mathcal{A}}]$. 
We shall show $\sigma_0 \not \in cl_b({\mathcal{A}})$.
Now 
$$cl_df[{\mathcal{A}}]=\cap \{D' \colon f[{\mathcal{A}}] \subseteq D'\}.$$
In particular, there is some $D_0 \in {\mathcal{D}}$ with
$f[{\mathcal{A}}] \subseteq D'_0$, yet $f(\sigma_0) \not \in D'_0$. 
Now 
$$D'_0 = \{\Pi_d(\gamma) \colon D_0 \in \Pi_d(\gamma)\}.$$
So $\sigma \in {\mathcal{A}}$ implies $f(\sigma) \in D'_0$, whence 
$D_0 \in f(\sigma)$ and $f^\leftarrow(D_0) \in \sigma$. 
Therefore $\sigma \in [f^\leftarrow (D_0)]'$. 
Hence ${\mathcal{A}} \subseteq [f^\leftarrow (D_0)]'$, closed in $T_bX$,
so $cl_b({\mathcal{A}}) \subseteq [f^\leftarrow(D_0)]'$.
But $f(\sigma_0) \not \in D'_0$. 
Therefore $D_0 \not \in f(\sigma_0)$ and 
$\sigma_0 \not \in [f^\leftarrow (D_0)]'$.
\end{proof}

\nocite{*}
\providecommand{\bysame}{\leavevmode\hbox to3em{\hrulefill}\thinspace}
\providecommand{\MR}{\relax\ifhmode\unskip\space\fi MR }
\providecommand{\MRhref}[2]{%
  \href{http://www.ams.org/mathscinet-getitem?mr=#1}{#2}
}
\providecommand{\href}[2]{#2}

\end{document}